\documentclass[10pt,a4paper]{amsart}

\usepackage[foot]{amsaddr}
\usepackage{amsfonts,epsfig,stmaryrd}
\usepackage{newlfont,amsfonts,amssymb,amsmath,amsthm,amsgen,amscd,datetime,dsfont,hyperref,tensor}

\usepackage[small,nohug,heads=littlevee]{diagrams}
\usepackage[normalem]{ulem}
\diagramstyle[labelstyle=\scriptstyle]

\usepackage{multicol,picinpar,enumerate,cleveref,mathabx}
\usepackage{euscript,color}

\newcommand{\hol}{\mathfrak{hol}}

\newcommand{\+}{\oplus}

\newcommand{\so}{\mathfrak{so}}

\newcommand{\hook}{\makebox[7pt]{\rule{6pt}{.3pt}\rule{.3pt}{5pt}}\,}

\renewcommand{\d}{\mathrm{d}}

\newcommand{\End}{\mathrm{End}}

\newcommand{\R}{\mathbb{R}} 

\renewcommand{\L}{\mathcal L}

\renewcommand{\k}{\mathfrak{kill}}

\theoremstyle{definition}
\newtheorem{definition}{Definition}

\newtheorem*{example*}{Example}
\newtheorem*{remark*}{Remark}

\theoremstyle{plain}

\newtheorem*{lem*}{Lemma}

\newtheorem{theorem}[definition]{Theorem}
\newtheorem*{theorem*}{Theorem}

\numberwithin{equation}{section}

\begin{document}

\title{Killing vector fields on semi-Riemannian product manifolds}

\thanks{This work was supported by 
 the Australian Research
Council (Discovery Program DP190102360).
 }
 \author[Federico Costanza]{Federico Costanza}
\email{efcostanza@gmail.com}

\author[Thomas Leistner]{Thomas Leistner}\address[Thomas Leistner, corresponding author]{School of Mathematical Sciences, University of Adelaide, SA~5005, Australia}\email{thomas.leistner@adelaide.edu.au}

\dedicatory{\small Dedicated to the memory of Joseph A.~Wolf
}


\begin{abstract}
Hano's theorem states that the space of Killing vector fields of a complete simply connected Riemannian manifold is isomorphic to the direct sum of the Killing vector fields of the factors in its de Rham decomposition. We prove a  generalisation of this theorem to manifolds with indefinite metrics that requires an assumption on  the factors, and show by example  why this assumption is needed.
 \end{abstract}

\maketitle
\setcounter{tocdepth}{1}

\section{Introduction}
If a semi-Riemannian manifold $(M,g)$ is isometric to a product of semi-Riemannian manifolds, many of its geometric properties can be analysed by studying the geometry of the factors. One of the geometric properties that can be analysed in this way is the existence of infinitesimal symmetries, or more precisely, of Killing vector fields. A classical theorem by Hano \cite{Hano55} states that for a simply connected and complete Riemannian manifold $(M,g)$  that is decomposed according to its de Rham decomposition \cite{derham52}, the space of Killing vector fields of $(M,g)$ is given by the sum of the Killing vector fields of the factors. In the present article we will show to which extent Hano's theorem can be generalised to semi-Riemanniann manifolds, i.e.~manifolds where the metric is not necessarily {\em positive} definite. Our result in Theorem \ref{thm} requires an assumption that prevents the existence of parallel vector fields on (one of) the factors. We will also show by an example that this assumption is indeed needed. The proof  of Theorem \ref{thm} uses the notion of the  {\em Killing connection}, which arises from the natural prolongation of the Killing equation. It  was first considered  in \cite{Kostant55}, and was recently studied systematically  in \cite{CostanzaPhD,CostanzaEastwoodLeistner21,calabi2}, where it turned out to be a useful tool to study the Killing operator. 
Our proof is based on the integrability condition for the existence of a Killing vector field that arises from the curvature of the Killing connection.


\subsubsection*{Acknowledgements.}
The result of this paper is an extension of Theorem 2.1.6 in the first author's PhD thesis \cite{CostanzaPhD} which was written under supervision of the second author. We would like to thank Michael Eastwood for taking the role of co-supervisor, Ben McMillan for participating in  weekly meetings, and both for helpful discussions and invaluable inspiration. 

We would like to thank the organisers of the 
workshop {\em Spectrum and Symmetry for Group Actions in Differential Geometry II}  and everyone at the MATRIX Mathematics Research Institute in Creswick for their hospitality. During  the two weeks of this workshop
we  enjoyed Joe Wolf’s company, his wisdom, and his wit, shortly before his passing.  At the workshop, Joe gave the last talk, and, as always, he inspired everyone.

\section{The de~Rham--Wu decompositon}
The fundamental result about the splitting of a semi-Riemannian manifold as a product according to its holonomy representation is the decomposition theorem of de~Rham and Wu \cite{derham52,wu64}. In order to state it we need to recall some notions. 
A semi-Riemannian manifold $(M,g)$ is {\em indecomposable} if it is not isometric to a semi-Riemannian product. For complete and simply connected manifolds this is equivalent to the property that the holonomy group of $(M,g)$ at $p\in M$  acts {\em indecomposably} on a tangent space $T_pM$, i.e.~without a non-degenerate invariant subspace.  A simply connected and complete {\em Riemannian} manifold, i.e.~when $g$ is positive definite, is indecomposable if and only if the holonomy group acts irreducibly,  i.e.~without any non-trivial invariant subspace at all. This holds since positive definite metrics do not allow for subspaces of the tangent space on which the metric degenerates.

The de~Rham--Wu decomposition theorem now states the following: a semi-Riemannian manifold $(M,g)$ that is geodesically complete and simply connected is isometric to a semi-Riemannian product $(M_0,g_0)\times (M_1,g_1)\times \ldots \times (M_k,g_k)$ of simply connected, complete semi-Riemannian manifolds, where $(M_0,g_0)$ is flat (or empty) and the $(M_i,g_i)$ are indecomposable for $i=1, \dots, k$ (or empty).
For Riemannian manifolds, i.e.~when $g$ is positive definite and hence the $(M_i,g_i)$ are irreducible for $i\ge 1$, this theorem was proven by G.~de~Rham \cite{derham52}. Hence, when applied to Riemannian manifolds, we just refer to it as the {\em de~Rham decomposition}. For semi-Riemannian manifolds, i.e.~when $g$ can be indefinite, the theorem was proven by 
H.~Wu
 \cite{wu64}, and in the general case we refer to its result as the {\em de~Rham--Wu decomposition}.

\section{Hano's theorem and its generalisation to indefinite metrics}
 A  Killing vector field of a semi-Riemannian manifold $(M,g)$ is a vector field $\xi$ whose flow preserves the metric $g$, or equivalently, which satisfies  $\L_\xi g=0$. Here $\mathcal L_\xi$ denotes the Lie derivative of a tensor along $\xi$, i.e.~
\[(\L_\xi g)(X,Y)=\xi(g(X,Y))-g([\xi,X],Y)-g(X,[\xi,Y])=g(\nabla_X\xi,Y)+g(X,\nabla_Y\xi),\]
where $\nabla$ is the Levi-Civita connection of $g$. With respect to the Lie bracket of vector fields, the Killing vector fields of $(M,g)$ form a real Lie algebra $\mathfrak{kill} (M,g)$. Note that a vector field $\xi$ is Killing if and only if its  dual one-form  $\xi^\flat:=g(\xi,.)$ satisfies
\[ \nabla \xi^\flat = \d\xi^\flat,\]
i.e.~if and only if $\nabla \xi^\flat$ is a $2$-form.

For Riemannian manifolds the relation between Killing vector fields and the de~Rham decomposition is well-known due to a quite satisfactory theorem by Hano.
\begin{theorem*}[Hano \cite{Hano55}]
Let $(M,g)=(M_0,g_0)\times (M_1,g_1)\times \ldots \times (M_k,g_k)$ be the de~Rham decomposition of a complete simply connected Riemannian manifold. Then the Lie algebra of Killing vector fields decomposes as
\[ \k(M,g)= \k(M_0,g_0)\+\k (M_1,g_1)\+ \ldots \+ \k(M_k,g_k).\]
\end{theorem*}
In  \cite{Hano55} and also in \cite[Theorem VI.3.5]{ko-no1} this theorem  is stated in terms of connected components of isometry groups, however we formulate it here on the Lie algebra level. Moreover, in both references the same splitting holds for affine vector fields of a Riemannian manifold.

Note that this theorem does not imply that the Killing vector fields of a product manifold split as direct sum of the Killing vector fields of the factors: clearly Euclidean space of dimension $n>1$ decomposes as a product, however its Killing vector fields do not split because of the presence of infinitesimal rotations. This is the reason why in Hano's theorem the flat factor is not split further into (irreducible) lines. We will come back to this feature.

Although it has been claimed in the literature, for example in \cite[Theorem 3]{BatatCastrillon-LopeRosado-Maria15}, Hano's theorem does not directly generalise to all semi-Riemannian manifolds and their de~Rham--Wu decomposition. This has been noted for example in \cite[Remark 3.6]{LeistnerTeisseire22}, where it is observed that a product of a Euclidean space and a  Cahen--Wallach space has more Killing vector fields than the sum of the Killing vector fields of the factors. A {\em Cahen--Wallach space} is a simply connected indecomposable, non-irreducible Lorentzian symmetric space. It is given by the following metric on $\R^{n+2}\ni (t,v,x^1,\ldots, x^n)$,
\[g=2 \d t ( \d v  +  x^iQ_{ij}x^j \d t)+ \delta_{ij}\d x^i \d x^j,\]
where  $(Q_{ij})$ is a non-degenerate symmetric $n\times n$-matrix. A Cahen--Wallach space $(\R^{n+2},g)$ is indecomposable, but admits the parallel null (i.e.~lightlike, or isotropic) vector field $\partial_v$, and hence is not irreducible.

Given these examples, one is lead to believe that the Euclidean factor  causes a generalisation of Hano's theorem to fail. However, the following example shows that it is not that simple.

\begin{example*}
Let $(M_+,g_+)$ and $(M_-,g_-)$ be two indecomposable semi-Riemannian manifolds that both admit a parallel null vector field $\xi_+$ and $\xi_-$, respectively. For example, this could be two Cahen--Wallach spaces. Since both $\xi_\pm$ are parallel, their metric duals $\xi_\pm^\flat\in \Gamma(T^*M_\pm)$ are closed. Hence, locally, or globally if we assume $M_+$ and $M_-$ to be simply connected, there are two functions $f_\pm\in C^\infty(M_\pm)$ such that $\d f_\pm=\xi_\pm^\flat$.

Now consider the semi-Riemannian product manifold $(M,g)=(M_+,g_+)\times (M_-,g_-)$. We can pull back  $\xi_\pm$, $\xi_\pm^\flat$ and $f_\pm$ from the factors to $(M,g)$, and, by abusing notation in a standard way, we denote them by the same symbols. Then the vector fields $\xi_\pm$'s are parallel  on $(M,g)$, and so are the $1$-forms $\xi_\pm^\flat$. Consider the vector field
\begin{equation}
\label{exampleVF}\xi =f_+  \xi_- -f_-\xi_+\end{equation}
on $(M,g)$, with covariant derivative
\[\nabla \xi^\flat = \d f_+ \otimes \xi_-^\flat -\d f_-\otimes \xi_+^\flat = \xi_+\wedge \xi_-.\]
 Hence, $\xi$ is a Killing vector field on $(M,g)$. 
 
 If we now assume both $(M_i,g_i)$ to be simply connected and geodesically complete, then $(M_+,g_+)\times (M_-,g_-)$ is the de~Rham--Wu decomposition of $(M,g)$, however $\xi$  is not a sum of Killing vector fields on $M_+$ and $M_-$. This example shows that a {\em direct} generalisation of Hano's theorem to indefinite semi-Riemannian manifolds cannot hold, not even when a flat factor is excluded from the de~Rham--Wu decomposition, i.e.~when only indecomposable factors are allowed.
\end{example*}

%
%
%

With the lesson learned from this example, in our generalisation of Hano's theorem we cannot allow the existence of a parallel vector field on the factors.
Our result is the following.
\begin{theorem}
\label{thm}
Let $(M,g)=(M_+,g_+)\times (M_-,g_-)$ be a product of two semi-Riemannian manifolds and such that $(M_+,g_+)$ is real analytic, simply connected, and does not admit a parallel vector field.
Then the Lie algebra of Killing vector fields decomposes as
\[ \k(M,g)= \k(M_+,g_+)\+\k (M_-,g_-).\]
In particular, if $(M,g)= (M_1,g_1)\times \ldots \times (M_k,g_k)$ is the de~Rham--Wu decomposition  of a complete simply connected semi-Riemannian manifold such that all $(M_i,g_i)$ are real analytic and at most one factor admits a parallel vector field, then the Lie algebra of Killing vector fields decomposes as
\[ \k(M,g)= \k (M_1,g_1)\+ \ldots \+ \k(M_k,g_k).\]
\end{theorem}
A version of this result can be found  in the first author's PhD thesis \cite[Theorem 2.1.6]{CostanzaPhD}. There,  real analyticity is not required, but it is assumed  that the curvature tensor of $(M_+,g_+)$ has {\em zero nullity}, i.e. that at each $p\in M_+$ it holds that $\{X\in T_pM_+\mid X\hook R_+=0\}=\{0\}$. This assumption is stronger  than the non-existence of a parallel vector field. 
Of course,  in an optimal version of Theorem~\ref{thm} the assumption of real analyticity should be dropped, and we believe that it is not needed. 
 This however involves a more detailed study of the involved holonomy algebras and hence is work in progress. Unfortunately, the assumption of real analyticity  prevents us from obtaining Hano's theorem as a corollary to Theorem~\ref{thm}.

Note  that for indefinite metrics  the factors without parallel vector fields in the de~Rham--Wu decompositon do not have to be irreducible. A general manifold  of Walker type (see for example \cite{walkerbook}) will have holonomy invariant null planes without having any parallel null vector fields. For example, an indecomposable Lorentzian manifold with a recurrent vector field $\xi$, i.e.~a null vector field with
\[\nabla \xi= \theta\otimes \xi,\quad \text{ where $\theta$ is  a {\em non-closed} $1$-form,}\]
is not irreducible, but also does not admit a parallel null vector field. 
Theorem \ref{thm} applies to manifolds that have such a factor in their de~Rham--Wu decomposition, and therefore it is not just a starightforward analogue of Hano's theorem for manifolds that decompose into irreducible factors. 

Our proof of Theorem~\ref{thm} in the last section uses  the {\em Killing connection},  that arises from the prolongation of the Killing equation. We will introduce the details in the next section.

 \section{The Killing connection}
 
In order to prove Theorem \ref{thm} we will present a vector bundle with connection whose parallel sections are in one to one correspondence with Killing vector fields. This was  originally introduced by Kostant in \cite{Kostant55} and later studied by Geroch \cite{Geroch69} in the context of conformal Killing transport in general relativity. More recently, this connection was systematically studied in \cite{CostanzaPhD,CostanzaEastwoodLeistner21,calabi2}, used to analyse the range of the Killing operator, and the term {\em Killing connection} was coined to denote it. 

Let $(M,g)$ be a semi-Riemannian manifold with Levi-Civita connection $\nabla $ and 
 curvature tensor
\[R(X,Y)=\left[\nabla_X,\nabla_Y\right]-\nabla_{[X,Y]}.\]
Within the endomorphism bundle $\End(TM)$, we define the vector bundle
\[\so(TM, g):=\{A\in \End (TM)\mid A\cdot g=0\},\]
where $A\cdot g$ denotes the  natural action of $A$ on the tensor field $g$,
\[(A\cdot g)(X,Y)=-g(AX,Y)-g(X,AY).\]
We will call  $E : = TM \oplus \so(TM, g)$  the \textit{Killing bundle} of $(M, g)$, and $D : \Gamma(E) \to \Gamma(\wedge^{1}M \otimes E)$, defined by
\[
D_{X}
\begin{pmatrix}
	\xi \\
	A
\end{pmatrix}
:=
\begin{pmatrix}
\nabla_X \xi + AX \\
\nabla_{X}A + R(X, \xi)
\end{pmatrix},
\]
will be the \textit{Killing connection} on $E$. For a given vector field $\xi$ of $M$, we define the endomorphism $A_{\xi} : = - \nabla \xi$, of $TM$, which coincides with $\mathcal{L}_{\xi} - \nabla_{\xi}$ when acting on vector fields, where $\cal L_\xi $ denotes the Lie derivative along $\xi$. Therefore, $\xi$ will be a Killing vector field of $(M, g)$ if and only if $A_{\xi}$ is in 
$\so(TM, g)$, since
\[
\mathcal{L}_{\xi} g = \nabla_{\xi}g + A_{\xi} \cdot g = A_{\xi} \cdot g.
\]
A straightforward computation, that was first carried out in \cite{Kostant55}, shows that if $\xi$ is a Killing vector field, then
\begin{equation}\label{nabA}
\nabla_{X} A_{\xi} = - R(X, \xi) \quad \text{for all} \quad X \in \Gamma(TM),
\end{equation}
and that, consequently, the parallel sections of $E$ are exactly those of the form $(\xi, A_{\xi})$ such that $\xi$ is a Killing vector field of $(M, g)$. 
Hence, the vector space of  parallel sections of $(E,D)$ is isomorphic to $\k(M,g)$ under the isomorphism 
\begin{equation}\label{eqn:Killing_isomorphism}
 \xi \to \phi_{\xi}:= (\xi, A_{\xi})= (\xi, -\nabla {\xi}).
\end{equation}
The curvature of the Killing connection is the homomorphism of vector bundles $\kappa : E \to \Lambda^{2}M \otimes E$, given by
\begin{equation}\label{eqn:Killing_curvature}
\kappa(X, Y)
\begin{pmatrix}
\xi \\
A
\end{pmatrix}
=
-
\begin{pmatrix}
0 \\
(\nabla_{\xi}R)(X, Y)+ (A \cdot R)(X, Y)
\end{pmatrix}.
\end{equation}
Now we have all the ingredients for the proof of Theorem~\ref{thm}.

\section{The proof of Theorem \ref{thm}}
Let $(M, g) = (M_+, g_+) \times (M_-, g_-)$ be a product of semi-Riemannian manifolds with curvature tensors $R_+$ and $R_-$.
We denote the pull-backs of $g_\pm$ and $R_\pm$ by the natural projection $\pi_{\pm} : M_{+} \times M_{-} \to M_\pm$ to $M_\pm$ by the same symbols.

It is clear that a the pull-back to $M$ of a Killing vector field on $M_+$ or $M_-$ is a Killing vector field of $(M,g)$. Hence $\k(M_+,g_+)\+\k(M_-,g_-)$ injects into $\k(M,g)$. We have to show that any Killing vector field of $(M,g)$ is of this form, provided that
$(M_+, g_+)$ satisfies the assumptions of the theorem. We will use the Killing connection for this.

On  $(M, g) = (M_+, g_+) \times (M_-, g_-)$ the Killing bundle will decompose as
\[
E =
\begin{matrix}
TM_+ & & TM_- & & \\
\oplus & \oplus & \oplus & \oplus & TM_+ \wedge TM_-\\
\mathfrak{so}(TM_+, g_+) & & \mathfrak{so}(TM_-, g_-) & & 
\end{matrix},
\]
where we identify $TM_+\wedge TM_-$ with a subbundle of  $\so(TM,g)$ via
\[ v_+\wedge v_-\longmapsto g(v_+,.)\otimes v_- - g(v_-,.)\otimes v_+.\]
Again we have abused notation and have identified $TM_\pm \oplus \so(TM_\pm, g_\pm)$ with the pullback of the Killing bundle of $(M_\pm, g_\pm)$ with respect to the natural projection $\pi_{\pm} : M_{+} \times M_{-} \to M_\pm$. It is evident that both $E_\pm:=TM_\pm\+\so(TM_\pm,g_\pm)$ are parallel subbundles of $E$ whereas $TM_+ \wedge TM_-$ is not.

Let $\xi$ be a Killing vector field of $(M,g)$ and denote by $(\xi, A)\in \Gamma(E)=\Gamma(TM\+\so(TM,g))$ 
 the corresponding parallel section with $A=-\nabla \xi$.  We will show that $A$ is a section of $\mathfrak{so}(TM_+, g_+) \+\mathfrak{so}(TM_-, g_-) $, i.e.~ that it has no $TM_+\wedge TM_-$-component.
Since $(\xi, A)$ is parallel, we have
\[
0=\kappa(X,Y)\left(\xi \atop A\right)\ = - \ 
\left(\begin{array}{c}0\\ (\nabla_\xi R)(X,Y) + (A\cdot R)(X,Y)\end{array} \right),\]
so that, when writing out the action $A\cdot R$ of $A$ on the curvature tensor $R$, 
\begin{equation}
\label{killcurv} (\nabla_\xi R)(X,Y) + [A, R(X,Y)]- R(AX,Y)-R(X,AY)=0.\end{equation}
Now we use the fundamental property of the curvature tensor (and its derivatives) of a product manifold. For each $k=0,1,\ldots$ we have
\begin{equation}\label{prodR}
(\nabla_{Z_k} \ldots \nabla_{Z_1}R)(X,Y)Z=\left\{\begin{array}{ll}
(\nabla_{Z_k} \ldots \nabla_{Z_1}R_+)(X,Y)Z,&\text{ if }X,Y,Z,Z_1, \ldots, {Z_k}\in TM_+,\\
(\nabla_{Z_k} \ldots \nabla_{Z_1}R_-)(X,Y)Z,&\text{ if }X,Y,Z,Z_1, \ldots, {Z_k}\in TM_-,\\
0,&\text{ otherwise.}
\end{array}\right\}.
\end{equation}
Then, when setting $X=X_+\in TM_+$ and $Y=X_-\in TM_-$ in equation~(\ref{killcurv}), we get
\[
0=R(AX_+,X_-)+R(X_+,AX_-)=R_-(AX_+,X_-)+R_+(X_+,AX_-).
\]
which yields the two equations
\begin{equation}
\label{mixed}
R_-(AX_+,X_-)=0,\qquad R_+(X_+,AX_-)=0, \quad\text{ for all $X_+\in TM_+$ and $X_-\in TM_-$.}
\end{equation}
We will now show inductively that the same holds for the $k$-th covariant derivatives,
\begin{equation}
\label{kmixed}
(\nabla^{(k)}R_-)(AX_+,X_-)=0,\quad (\nabla^{(k)} R_+)(X_+,AX_-)=0,\quad \text{ for all $k\ge 0$ and all $X_\pm\in TM_\pm$. }
\end{equation}
Clearly the base case $k=0$ is provided by equations~(\ref{mixed}), so we now assume that ~(\ref{mixed}) holds up to $k-1$ for $k\ge 1$.
Because of~(\ref{prodR}) it is sufficient to show that for each $X_+\in TM_+$ it is
\[
(\nabla_{Z_k} \ldots \nabla_{Z_1}R_-)(AX_+,X_-)=0,\quad \text{ for $X_-,Z_i\in TM_-$,}
\]
and for each $X_-\in TM_-$,
\[
(\nabla_{Z_k} \ldots \nabla_{Z_1}R_+)(AX_-,X_+)=0,\quad \text{ for $X_+,Z_i\in TM_+$.}
\]
We prove the second equation, and the first will follow, {\em mutatis mutandis}. We fix an $X_-\in TM_-$. Because of the induction hypothesis and since $TM_+$ is a parallel distribution, for $Z_i\in TM_+$,  it is
\[
(\nabla_{Z_k} \ldots \nabla_{Z_1}R_+)(AX_-,X_+)
=
-
(\nabla_{Z_{k-1}} \ldots \nabla_{Z_1}R_+)(\nabla_{Z_k}(AX_-),X_+).
\]
This we compute, using  $A=-\nabla\xi$ and the resulting equation~(\ref{nabA}), as
\[
\nabla_{Z_k}(AX_-)= (\nabla_{Z_k}A)X_- +A\nabla_{Z_k}X_-=-R(Z_k,\xi)X_-=0,\]
because of~(\ref{prodR}), $Z_k\in TM_+$, and since we can choose $X_-$ such that $\nabla X_-|_{TM_+}=0$. This proves equations~(\ref{kmixed}) for all $k\ge 0$.

The proof concludes as follows. Equations~(\ref{kmixed}) imply that on both factors we have that 
\begin{equation}
\label{kmixed+-}
(\nabla_-^{(k)}R_-)(AX_+,X_-)=0,\qquad (\nabla_+^{(k)} R_+)(X_+,AX_-)=0,\quad \text{ for all $k\ge 0$, }
\end{equation}
where $\nabla_\pm$ are the Levi-Civita connections of $g_\pm$.
Under our assumptions  on $(M_+,g_+)$
 we would like to conclude from this that $\mathrm{pr}_{TM_+}(AX_-)=0$. We prove this by contradiction. Assume that there is a $p\in M_+$ and $0\not=Y_+\in T_pM_+$ such that 
 \[(\nabla_+^{(k)} R_+)(X_+,Y_+)|_p=0,\quad \text{ for all $X_+\in T_pM_+$ and all $k\ge 0$.}\]
By the symmetries of $R_+$ and since $(M_+,g_+)$ is assumed to be real analytic, this implies that the holonomy algebra of $(M_+,g_+)$ at $p$ annihilates $Y_+$,
 \[\hol_p(M_+,g_+)\cdot Y_+=0,\]
see for example \cite[Section II.10]{ko-no1} for the equality between  the infinitesimal and full holonomy algebra in the real analytic case.
Since $M_+$ is assumed to be simply connected, it implies that $Y_+$ can be extended to all of $M_+$ to a parallel vector field, which contradicts the assumption. As a result, we get the desired $\mathrm{pr}_{TM_+}(AX_-)=0$.
This implies 
 $\mathrm{pr}_{TM_+} \circ A|_{TM_-}=0$, but also that  $\mathrm{pr}_{TM_-} \circ A|_{TM_+}=0$, so that $A$ is a section of $\so(TM_+,g_+)\+\so(TM_-,g_-)\subset \so(TM,g)$. 
Hence, \[\left(\xi, A\right)=\left(\xi_+, A_+\right)+\left(\xi_-, A_-\right)\ \in\  \Gamma(E_+\+E_-),\] where $E_\pm=TM_\pm\+\so(TM_\pm,g_\pm)$. Since the $E_\pm$'s are parallel subbundles in $E$,
both $\left(\xi_\pm, A_\pm\right)$ are  parallel sections of $E_\pm$ and hence both $\xi_\pm$'s are Killing for $g_\pm$.
Therefore, with the isomorphism provided in equation (\ref{eqn:Killing_isomorphism}), we obtain the desired result
	\[
	\mathfrak{kill}(M, g) = \mathfrak{kill}(M_+, g_+) \oplus \mathfrak{kill}(M_-, g_-).
	\]
	The remainder of the theorem follows from repeating this process for each de Rham--Wu factor of $(M, g)$.

\bibliographystyle{abbrv}
\bibliography{GEOBIB}
\end{document}